\documentclass[11pt]{article}
\usepackage{amssymb,amsthm,amsmath,palatino}
\usepackage[pdftex]{graphicx}
\usepackage[letterpaper,hmargin=1in,vmargin=1.25in]{geometry}
\usepackage{url}

 \newtheorem{thm}{Theorem}[section]
 
 \newtheorem{lem}[thm]{Lemma}
 
 \newtheorem{conj}[thm]{Conjecture}
\theoremstyle{definition}
 \newtheorem{defn}[thm]{Definition}
 \theoremstyle{remark}

 \numberwithin{equation}{section}

\DeclareMathOperator{\rank}{rank}
\newcommand{\RR}{\mathbb{R}}
\newcommand{\Fq}{\mathbb{F}_q}
\newcommand{\uu}{\mathbf{u}}

\newcommand{\ww}{\mathbf{w}}
\newcommand{\xx}{\mathbf{x}}
\renewcommand{\d}{\mathrm{d}}

\title{Metric dimension of dual polar graphs}
\author{Robert~F.~Bailey\footnote{School of Science and the Environment (Mathematics), Grenfell Campus, Memorial University of Newfoundland, Corner Brook, NL A2H~6P9, Canada. E-mail: \texttt{rbailey@grenfell.mun.ca}} \and 
Pablo Spiga\footnote{Dipartimento di Matematica e Applicazioni, University of Milano-Bicocca, Via Cozzi 55, 20125 Milano, Italy. \mbox{E-mail:} \texttt{pablo.spiga@unimib.it}}}
\begin{document}

\maketitle

\begin{abstract}
A {\em resolving set} for a graph $\Gamma$ is a collection of vertices $S$, chosen so that for each vertex $v$, the list of distances from $v$ to the members of $S$ uniquely specifies $v$. The {\em \mbox{metric} dimension} $\mu(\Gamma)$ is the smallest size of a resolving set for $\Gamma$.  We consider the metric dimension of the {\em\mbox{dual polar graphs}}, and show that it is at most the rank over $\RR$ of the incidence matrix of the corresponding polar space. We then compute this rank to give an explicit upper bound on the metric dimension of dual polar graphs, as well as the halved dual polar graphs.
\end{abstract}

\section{Introduction} \label{section:intro}

\subsection{Resolving sets and metric dimension} \label{section:intro-metdim}

Let $\Gamma$ denote a graph with vertex set $V$ and edge set $E$, which we assume to be finite, connected, loopless, and with no multiple edges.  A {\em resolving set} for $\Gamma$ is a subset $S\subseteq V$ with the property that, for every $u\in V$, the list of distances from $u$ to each of the elements of $S$ uniquely identifies $u$; equivalently, for two distinct vertices $u,w\in V$, there exists $x\in S$ for which $\d(u,x)\neq\d(w,x)$. (Here, $\d(x,y)$ denotes the length of a shortest path from $x$ to $y$ in $\Gamma$.) The {\em metric dimension} of $\Gamma$ is the smallest size of a resolving set for $\Gamma$, and we denote this by $\mu(\Gamma)$.  These notions were introduced to graph theory in the 1970s by Slater~\cite{Slater75} and, independently, by Harary and Melter~\cite{Harary76}; in more general metric spaces, the concept can be found in the literature much earlier (see~\cite{Blumenthal53}).  

The 2011 survey article by Cameron and the first author~\cite{bsmd} was the first to consider, systematically, the metric dimension of distance-regular graphs; in that paper it was observed that the work of Babai from 1981~\cite{Babai81} gives a general upper bound of $O(\sqrt{n} \log n)$ on the metric dimension of a primitive distance-regular graph on $n$ vertices (see~\cite[Theorem 3.15]{bsmd}).  Since then various papers have been written on this subject, as well as the more general problem of class dimension of association schemes, often focused on particular families: for example, Johnson and Kneser graphs~\cite{jk}, Grassmann graphs~\cite{grassmann}, bilinear forms graphs~\cite{FengWang}, symplectic dual polar graphs~\cite{symplectic}, incidence graphs of designs and geometries~\cite{imprimitive,incidence,Bartoli,PGnq,PG-lines,Heger2,Heger}, and point graphs of generalized quadrangles~\cite{Gravier15,HegerGQ}; see also~\cite{small,singular,attenuated,flat,GuoWangLi,fourfamilies}.  In this paper, we will consider {\em dual polar graphs}, which are defined below.

\subsection{Polar spaces and dual polar graphs}

Let $V(n,q)$ denote the vector space of dimension $n$ over $\Fq$, the finite field with $q$ elements, equipped with either a symplectic, quadratic or Hermitian non-degenerate structure. Recall that to define one of these three structures on $V(n,q)$, we require a sesquilinear form $b:V(n,q)\times V(n,q)\to \mathbb{F}_q$ in two variables and/or a quadratic form $f:V(n,q)\to \mathbb{F}_q$ in one variable, where $f$ is either defined by $f(v)=b(v,v)$ for every $v\in V(n,q)$, or $b$ is obtained by polarizing $f$. A subspace $U$ of $V(n,q)$ is called {\em totally isotropic} if the restriction to $U$ of either the sequilinear form or the quadratic form is identically zero; that is, $b(v,w)=0$ or $f(v)=0$ for all $v,w\in U$. A {\em (classical) polar space} $\mathcal{P}$ is the collection of all totally isotropic subspaces of $V(n,q)$. The {\em Witt index} of $V(n,q)$ is the dimension of a largest totally isotropic subspace of $V(n,q)$.  The Witt index of $V(n,q)$ is often called the {\em rank} of $\mathcal{P}$.  The $1$-dimensional totally isotropic subspaces are the {\em points} of $\mathcal{P}$, while if $\mathcal{P}$ has Witt index $d$, the $d$-dimensional totally isotropic subspaces are the {\em generators} (or the {\em maximals}) of $\mathcal{P}$.

With this terminology, we make the following definition.

\begin{defn} \label{defn:polarspace}
Let $\mathcal{P}$ be a polar space over $\Fq$ with Witt index $d$.  The {\em dual polar graph} on $\mathcal{P}$ has the generators of $\mathcal{P}$ as vertices, and two generators are adjacent if and only if their intersection has dimension $d-1$.
\end{defn}

There are six families of classical polar spaces, and thus six families of dual polar graphs, arising from the classification of sequilinear and quadradic forms.  Numerical information about these spaces can be expressed in terms of the field order $q$, Witt index $d$, and a parameter $e\in\{0,\,1/2,\,1,\,3/2,\,2\}$ depending on the choice of form.  We note that a Hermitian polar space requires the field order $q$ to be a square.  These are summarized in the table below; our notation follows~\cite{BCN}.

\begin{table}[hbtp] \label{table:polar}
\centering
\renewcommand{\arraystretch}{1.25}
\begin{tabular}{c|c|c|c}
Polar space & Names &  Vector space & $e$ \\
\hline
$[C_d(q)] \cong Sp(2d,q)$ & Symplectic & $V(2d,q)$ & $1$ \\ 
$[B_d(q)] \cong \Omega(2d+1,q)$ & Orthogonal; parabolic quadric & $V(2d+1,q)$ & $1$ \\
$[D_d(q)] \cong \Omega^+(2d,q)$ & Orthogonal; hyperbolic quadric &  $V(2d,q)$ & $0$ \\
$[^2\! D_{d+1}(q)] \cong \Omega^-(2d+2,q)$ & Orthogonal; elliptic quadric &  $V(2d+2,q)$ & $2$ \\
$[^2\! A_{2d}(\sqrt{q})] \cong U(2d+1,q)$ & Unitary; Hermitian variety &  $V(2d+1,q)$ & $3/2$ \\
$[^2\! A_{2d-1}(\sqrt{q})] \cong U(2d,q)$ & Unitary; Hermitian variety &  $V(2d,q)$ & $1/2$
\end{tabular}
\vspace{1ex}
\renewcommand{\arraystretch}{1.0}
\caption{Classical polar spaces}
\end{table}

We refer the reader to Brouwer and van Maldeghem~\cite[Chapter~2]{BvM} or De Bruyn~\cite{DeBruyn} for background on polar spaces, and to Brouwer, Cohen and Neumaier~\cite[Section 9.4]{BCN} for dual polar graphs.  We will use the notation $\Gamma(q,d,e)$ to denote a dual polar graph when the type is unspecified.  The following result is taken from~\cite[Section~9.4]{BCN}.

\begin{lem} \label{lemma:dualpolarfacts}
Let $\Gamma(q,d,e)$ be the dual polar graph arising from a polar space $\mathcal{P}$.  Then the following holds.
\begin{itemize}
\item[(a)] The number of points of $\mathcal{P}$ is 
\[ \frac{(q^{d+e-1}+1)(q^d-1)}{q-1}. \]
\item[(b)] The number of generators of $\mathcal{P}$, and thus the number of vertices of $\Gamma(q,d,e)$, is
\[ \prod_{i=0}^{d-1} (q^{e+i}+1). \]
\item[(c)] If $U,W$ are vertices of $\Gamma(q,d,e)$, then $U,W$ are at distance $i$ if and only if $\dim(U\cap W)=d-i$.
\end{itemize}
\end{lem}
We note that when $e=1$, the dual polar graphs $[C_d(q)]$ and $[B_d(q)]$ have the same parameters, but are not isomorphic in general.

\section{The main theorem} \label{section:main}

In~\cite{jk,grassmann}, the first author, Meagher and others considered the metric dimension of Johnson graphs and Grassmann graphs, and obtained upper bounds on this equal to the rank of an appropriate incidence matrix.  This approach was subsequently used in~\cite{singular,attenuated,flat} for the class dimension of various families of association schemes.  In what follows, we shall adapt this technique for dual polar graphs.

Suppose that we have a dual polar graph $\Gamma(q,d,e)$ arising from a polar space $\mathcal{P}$.  For each $t\in\{1,\ldots,d\}$, let $\Omega_t$ denote the set of all totally isotropic $t$-dimensional subspaces in $\mathcal{P}$ (so that $\Omega_1$ is the set of points of $\mathcal{P}$, and $\Omega_d$ the set of generators of $\mathcal{P}$).

We recall (from~\cite{BvM}, for instance) that a graph is {\em strongly regular} with parameters $(n,k,a,c)$ if it has $n$ vertices, is regular with degree $k$, any pair of adjacent vertices have $a$ common neighbours, and any pair of non-adjacent vertices have $c$ common neighbours.  Polar spaces are a source of such graphs: the {\em collinearity graph} of a polar space $\mathcal{P}$ is the graph whose vertices are the points of $\mathcal{P}$, and where two distinct points are adjacent if and only if their span is totally isotropic.  The following facts will be of use to us.

\begin{lem} \label{lemma:collin}
Let $\Delta$ denote the collinearity graph of a polar space with parameters $q,d,e$ as above.  Then $\Delta$ is strongly regular with parameters $(n,k,a,c)$, where
\begin{align*}
n & = |\Omega_1| = (q^{d+e-1}+1)\frac{q^d-1}{q-1},\\
k & = q(q^{d+e-2}+1)\frac{q^{d-1}-1}{q-1}, \\
a & = (q-1) + q^2(q^{d+e-3}+1)\frac{q^{d-2}-1}{ q-1}, \\
c & = (q^{d+e-2}+1) \frac{q^{d-1}-1}{q-1}.
\end{align*}
Furthermore $\Delta$ has eigenvalues $\theta_0=k$, $\theta_1 = q^{d-1}-1$ and $\theta_2 = -(q^{d+e-2}+1)$, with multiplicities $m_0=1$, $m_1$ and $m_2$ (respectively), where
\begin{align*}
m_1 & = \frac{q^e(q^{d+e-2}+1)}{q^{e-1}+1}\frac{q^d-1}{q-1}, \\
m_2 & = \frac{q(q^{d+e-1}+1)}{q^{e-1}+1} \frac{q^{d-1}-1}{q-1}.
\end{align*}
\end{lem}

A proof is given by Brouwer and van Maldeghem~\cite[Theorem~2.2.12]{BvM}: the fact that $\Delta$ is strongly regular with the parameters $(n,k,a,c)$ as given is well-known, while the eigenvalues and their multiplicities can be calculated from the standard formulas for strongly regular graphs.

For a given $U\in\Omega_t$, the {\em incidence vector} of $U$ is the vector $\uu\in\RR^{\Omega_1}$ with entries $0$ or $1$ so that, for any $x\in\Omega_1$, the $x$-coordinate of $\uu$ is $1$ if $x\subseteq U$ and $0$ otherwise.  If we have a collection of totally isotropic subspaces $\mathcal{W}=\{W_1,\ldots,W_m\}$, the {\em incidence matrix} of $\mathcal{W}$ is the $m\times|\Omega_1|$ matrix whose rows are the incidence vectors of $W_1,\ldots,W_m$.  We shall refer to the incidence matrix of the collection of all generators of $\mathcal{P}$ as the {\em incidence matrix of $\mathcal{P}$}.  
The following lemma will be crucial in our upper bound on the size of a resolving set for $\Gamma(q,d,e)$.

\begin{lem} \label{lemma:rank}
Let $\mathcal{P}$ be a polar space with parameters $q,d,e$ as above. Then the incidence matrix of $\mathcal{P}$ has rank
\[ \frac{(q^{d+e-1}+1)(q^{d+e-1}-q^{e-1}+q-1)}{(q^{e-1}+1)(q-1)}. \]
\end{lem}

\proof Let $M$ be the incidence matrix of $\mathcal{P}$ and let $B=M^T M$. By standard linear algebra, $\rank(B)=\rank(M)$.  Clearly, $B$ is an $|\Omega_1| \times |\Omega_1|$ matrix with rows and columns indexed by the points of $\mathcal{P}$. For each $t\in \{0,\ldots,d-1\}$, let $N_t$ denote the number of generators containing a fixed $t$-dimensional subspace $U$. Observe that $U^\perp/U$ is a polar space of the same type as $\mathcal{P}$ with parameters $q,d-t,e$. Hence, using Lemma~\ref{lemma:dualpolarfacts}, we have that $N_t = \prod_{i=0}^{d-t-1} (q^{e+i}+1)$. Moreover, for two points $x,y\in \Omega_1$, the entry $B_{xy}$ is the number of generators containing both $x$ and $y$, so there are three possible values for this:
\[ B_{x,y} = \left\{ \begin{array}{ll}
   N_1 & \textnormal{if $x=y$}, \\
	 N_2 & \textnormal{if $x\neq y$ and the span of $x$ and $y$ is totally isotropic}, \\
	 0 & \textnormal{if $x\neq y$ and the span of $x$ and $y$ is not totally isotropic}.
   \end{array} \right. \]
From this, it follows that 
\[ B = N_1 I + N_2 A, \]
where $A$ is the adjacency matrix of the collinearity graph of $\mathcal{P}$.  Therefore, the eigenvalues of $B$ are $\lambda_i = N_1 + \theta_i N_2$, where $\theta_i$ is an eigenvalue of $A$, and $\lambda_i$ and $\theta_i$ have the same multiplicity.  The eigenvalues $\theta_0,\theta_1,\theta_2$ of $A$ and their multiplicities were given in Lemma~\ref{lemma:collin}, from which we observe that $N_1 + \theta_2 N_2 =0$.  Consequently, $B$ is a singular matrix, and its nullity is equal to $m_2$ (the multiplicity of $\theta_2$).  Therefore, we have
\[ \rank(M) = \rank(B) = |\Omega_1| - m_2 = \frac{(q^{d+e-1}+1)(q^{d+e-1}-q^{e-1}+q-1)}{(q^{e-1}+1)(q-1)}. \qedhere \] 
\endproof

By Lemma~\ref{lemma:dualpolarfacts}(c), a collection $\mathcal{S}=\{X_1,\ldots,X_m\}$ of generators of $\mathcal{P}$ will form a resolving set for $\Gamma(q,d,e)$ if and only if the map $\Omega_d \to \RR^m$ defined by $U \mapsto \left( \dim(X_1\cap U),\ldots,\dim(X_m\cap U) \right)$ is injective.  Now, for any totally isotropic subspaces $U,W$, we have that $\dim(U\cap W)=k$ if and only if $U$ and $W$ have $(q^k-1)/(q-1)$ points of $\mathcal{P}$ in common.  Thus we can phrase the ``resolving property'' for a collection of generators in terms of linear algebra: if $M$ is the incidence matrix of a collection of generators $\mathcal{S}=\{X_1,\ldots,X_m\}$, and $\uu$ is the incidence vector (as a column vector) of a given generator $U$, then the entry in position $i$ of the vector $M\uu \in\RR^m$ is $(q^{\dim(X_i\cap U)}-1)/(q-1)$.  Consequently, $\mathcal{S}=\{X_1,\ldots,X_m\}$ is a resolving set with incidence matrix $M$ if and only if, for any generators $U,W$ with incidence vectors $\uu,\ww$, we have that $M\uu=M\ww$ implies that $U=W$.

\begin{thm} \label{theorem:main}
Let $\Gamma(q,d,e)$ be the dual polar graph arising from a polar space $\mathcal{P}$.  Then the metric dimension of $\Gamma(q,d,e)$ is at most the rank, over $\RR$, of the incidence matrix of $\mathcal{P}$, that is, it is at most \[ \frac{(q^{d+e-1}+1)(q^{d+e-1}-q^{e-1}+q-1)}{(q^{e-1}+1)(q-1)}. \]
\end{thm}

\proof  Let $M$ denote the incidence matrix of $\mathcal{P}$.  Since $M$ is an $|\Omega_d| \times |\Omega_1|$ matrix, it has more rows than columns, and thus $\rank(M)\leq |\Omega_1|$.  By 
re-arranging rows if necessary, we can assume that
\[ M = \begin{pmatrix} A \\  B \end{pmatrix}, \]
where $A$ is a $\rank(M)\times |\Omega_1|$ matrix whose rows are linearly independent.  We will show that $A$ is the incidence matrix of a resolving set for $\Gamma(q,d,e)$.

By construction, we have $\rank(A)=\rank(M)$; since both matrices have the same number of columns, it follows that they have the same nullity.  However, if $M\xx=\mathbf{0}$ we must have $A\xx=\mathbf{0}$, and thus $\ker(M)\subseteq\ker(A)$; as $\ker(M)$ and $\ker(A)$ have the same dimension, it follows that $\ker(M)=\ker(A)$.  Now suppose that $U,W$ are generators of $\mathcal{P}$ with incidence vectors $\uu,\ww$, respectively.  Then we have
\begin{eqnarray*}
A\uu = A\ww & \iff & A(\uu-\ww) = \mathbf{0} \\ 
            & \iff & \uu-\ww \in \ker(A) \\
			& \iff & \uu-\ww \in \ker(M) \\
			& \iff & M\uu = M\ww.
\end{eqnarray*}
Therefore we must have $\dim(U\cap Z)=\dim(W\cap Z)$ for all $Z\in\Omega_d$.  In particular, this holds for $Z=U$, so $\dim(W\cap U)=\dim(U)$ and thus $U=W$.

Hence $A$ is indeed the incidence matrix of a resolving set for $\Gamma(q,d,e)$, of size $\rank(M)$.  The rest of the proof follows immediately from Lemma~\ref{lemma:rank}. \endproof

In Table~\ref{table:bounds} we restate the bound on the metric dimension according on the type of the polar space.

\begin{table}[hbtp] \label{table:bounds}
\centering
\renewcommand{\arraystretch}{2.25}
\begin{tabular}{c|c}
Graph  & Bound on metric dimension \\
\hline
$[C_d(q)]$, $[B_d(q)]$      & $\displaystyle{ \frac{1}{2} \frac{(q^d+1)(q^d+q-2)}{q-1} }$\\ 								
$[D_d(q)]$                  & $\displaystyle{ \frac{(q^{d-1}+1)(q^d+q^2-q-1)}{q^2-1} }$ \\
$[^2\! D_{d+1}(q)]$         & $\displaystyle{ \frac{q^{2(d+1)}-1}{q^2-1} }$ \\ 
$[^2\! A_{2d}(\sqrt{q})]$   & $\displaystyle{ \frac{(q^{d+\frac{1}{2}}+1)(q^{d+\frac{1}{2}}+q-\sqrt{q}-1)}{(\sqrt{q}+1)(q-1)} }$ \\ 
$[^2\! A_{2d-1}(\sqrt{q})]$ & $\displaystyle{ \frac{(q^{d-\frac{1}{2}}+1)(q^d+q^\frac{3}{2}-\sqrt{q}-1)}{(\sqrt{q}+1)(q-1)} }$
\end{tabular}
\vspace{1ex}
\renewcommand{\arraystretch}{1.0}
\caption{Bounds for each family of dual polar graphs}
\end{table}

\section{Halved dual polar graphs} \label{section:halved}

There are a number of other families of distance-regular graphs related to the dual polar graphs, as detailed in~\cite[Section 9.4C]{BCN}.  In this section, we adapt our results to one of these families, the halved dual polar graphs.  The graphs $[D_d(q)]$ arising from the hyperbolic quadric $\Omega^+(2d,q)$ are bipartite: this can be deduced from their parameters as distance-regular graphs, but can also be seen more directly from the geometry.

In $\Omega^+(2d,q)$, there is an equivalence relation $\sim$ on the set of all generators, where $U\sim W$ if and only if $d-\dim(U\cap W)$ is even.  This relation on generators has two equivalence classes, which are referred to as {\em Latins} and {\em Greeks}.  (Since the group $\mathrm{PSO}^+(2d,q)$ acts transitively on the generators, this labelling is arbitrary.)  On, say, the set of Latin generators, we may define the {\em halved dual polar graph} $\Delta(q,d)$, where $U$ and $W$ are adjacent if and only if they are at distance~$2$ in $\Gamma(q,d,0)$, i.e.\ if $\dim(U\cap W)=d-2$.  Distances between vertices in $\Delta(q,d)$ are precisely half those in $\Gamma(q,d,0)$, so two Latin generators $U$ and $W$ are at distance~$i$ in $\Delta(q,d)$ if and only if $\dim(U\cap W)=d-2i$.  We note that the geometries associated with the Latin or Greek generators are sometimes known as {\em half-spin geometries} (see~\cite{DeBruyn}, for instance).

Because there is a one-to-one correspondence between distances and dimensions of intersections, the methods of Section~\ref{section:main} may be easily adapted to bound the metric dimension of $\Delta(q,d)$, as follows.  Once again, the key observation is the rank of the corresponding incidence matrix.

\begin{lem} \label{lemma:latin-greek-rank}
Let $M$ be the incidence matrix of $\Omega^+(2d,q)$, and let $M_1$ and $M_2$ be the incidence matrices of the Latin and Greek generators in $\Omega^+(2d,q)$, respectively.  Then
\[ \rank (M_1) = \rank (M_2) = \rank (M) =  \frac{(q^{d-1}+1)(q^d+q^2-q-1)}{q^2-1}. \]
\end{lem}

\proof Since the choice of which were the Latin and Greek generators was arbitrary, it is clear that $\rank(M_1)=\rank(M_2)$.  Also, $\rank(M)$ was calculated in Lemma~\ref{lemma:rank}, so it remains to show that $\rank(M_1)=\rank(M)$.  

As in the proof of Lemma~\ref{lemma:rank}, we let $N_t$ denote the number of generators of $\Omega^+(2d,q)$ containing a fixed $t$-dimensional subspace $U$, and we now let $N'_t$ denote the number of Latin generators containing $U$.  The totally isotropic subspaces through $U$ form a hyperbolic quadric $\Omega^+(2(d-t),q)$, so the generators of $\Omega^+(2d,q)$ through $U$ are in one-to-one correspondence with the generators of this $\Omega^+(2(d-t),q)$.  Consequently, the numbers of Latin and Greek generators containing $U$ must be equal, so we must have $N'_t=\frac{1}{2}N_t$.

Next, let $\Lambda=M_1^T M_1$, so $\Lambda$ is a matrix with rows and columns indexed by the points of $\Omega^+(2d,q)$, and where the entries are
\[ \Lambda_{xy} = \left\{ \begin{array}{ll}
   N'_1 & \textnormal{if $x=y$}, \\
	 N'_2 & \textnormal{if $x\neq y$ and the span of $x$ and $y$ is totally isotropic}, \\
	 0 & \textnormal{if $x\neq y$ and the span of $x$ and $y$ is not totally isotropic}.
\end{array} \right. \]
Therefore, we have $\Lambda= N'_1 I + N'_2 A = \frac{1}{2}B$, where $A$ is the adjacency matrix of the collinearity graph, and $B=M^T M$ is the matrix from the proof of Lemma~\ref{lemma:rank}.  Thus 
\[ \rank(M_1)=\rank(\Lambda)=\rank(B)=\rank(M) \]
and the result follows from Lemma~\ref{lemma:rank}.  \endproof

With this established, we now have a bound on the metric dimension for this class of graphs.

\begin{thm} \label{theorem:halved}
The metric dimension of the halved dual polar graph $\Delta(q,d)$ is at most
\[ \frac{(q^{d-1}+1)(q^d+q^2-q-1)}{q^2-1}. \]
\end{thm}

\proof Identical to the proof of Theorem~\ref{theorem:main}. \endproof

As explained in~\cite[Section 9.4C]{BCN}, the graph $\Delta(q,d)$ may also be constructed from the dual polar graph on the parabolic quadric $[B_{d-1}(q)]$, by forming its distance-1 or 2 graph.  Since the symplectic dual polar graph $[C_{d-1}(q)]$ has the same parameters as $[B_{d-1}(q)]$, it follows that its distance-1 or 2 graph is also distance-regular with the same parameters as $\Delta(q,d)$: this is known as the {\em Ustimenko graph}.  Also, the Ustimenko graphs arise as the halved graphs of the {\em Hemmeter graphs}, which have the same parameters as $[D_d(q)]$.  However, neither the Ustimenko nor Hemmeter graphs have vertex sets where distances are in one-to-one correspondence with dimensions of intersections of generators in a polar space, so there is no way to directly apply our methods to these graphs.  However, it seems reasonable to predict that the same upper bounds on metric dimension will apply, and we make the following conjecture.

\begin{conj}
The Ustimenko and Hemmeter graphs arising from $[C_{d-1}(q)]$ have metric dimension at most
\[ \frac{(q^{d-1}+1)(q^d+q^2-q-1)}{q^2-1}. \]
\end{conj}

\section{Final remarks} \label{section:remarks}

Previously, the only case of dual polar graphs of arbitrary diameter for which there is an existing bound in the literature is the symplectic dual polar graph $[C_d(q)]$, which was considered by Guo, Wang and Li~\cite{symplectic} in 2013. In that paper, the authors gave an upper bound on its metric dimension of $(q^d+1)(q^d+q-2)/(q-1)$, which is exactly double our bound.  In addition, our method is independent of the type of the polar space.

The bounds we obtained in Theorems~\ref{theorem:main} and~\ref{theorem:halved} are most effective when $d$ is large, as the upper bound on the metric dimension of $\Gamma(q,d,e)$ grows at a much slower rate when compared with the number of vertices.  For example, when $e=1$ (i.e.\ $[C_d(q)]$ and $[B_d(q)]$) the number of vertices is $O(q^{d(d+1)/2})$, but our upper bound is $O(q^{2d-1})$; Babai's general upper bound on the metric dimension (cf.\ \cite[Theorem~3.15]{bsmd}) evaluates as $O(q^{d(d+1)/4} \log q)$ here.  Similar analyses may be performed for the other values of $e$.

Our method is less effective for small $d$: for example, when $d=2$ and $e=1$ (i.e.\ $[C_2(q)]$ and $[B_2(q)])$, both the number of vertices and the upper bound from Theorem~\ref{theorem:main} are $O(q^3)$.  However, in these situations better bounds are already known: as well as Babai's bound, there are tighter bounds which make use of the geometry of the corresponding polar spaces.  Other than $[D_2(q)]$ (which is a complete bipartite graph), the dual polar graphs with $d=2$ are the line graphs of the classical generalized quadrangles.  By considering (i) the point graphs of the dual generalized quadrangles, and (ii) objects related to resolving sets known as {\em identifying codes}, bounds on the metric dimension of these graphs were found by Gravier {\em et al.}~\cite{Gravier15} and more recently by H\'eger {\em et al.}~\cite{HegerGQ}.  Notation for generalized quadrangles follows that used by Payne and Thas~\cite{PayneThas}.

\begin{itemize}
\item For $[B_2(q)]$, we have the lines of the parabolic quadric $Q(4,q)$, which is the dual of the symplectic GQ $W(q)$.  The point graph of $W(q)$ was considered in~\cite{HegerGQ}, with an upper bound of $4q$ on the metric dimension.  This compares with our upper bound of $\frac{1}{2}(q+2)(q^2+1)$. 
\item For $[C_2(q)]$, we have the lines of $W(q)$, which is the dual of $Q(4,q)$.  The point graph of $Q(4,q)$ was considered in~\cite{Gravier15}, with an upper bound of $5q-2$ on the metric dimension; for even $q$, this was improved in~\cite{HegerGQ} to $4q$.  In comparison, our upper bound is $\frac{1}{2}(q+2)(q^2+1)$. 
\item For $[^2\! D_3(q)]$, we have the lines of the elliptic quadric $Q^-(5,q)$, which is the dual of the Hermitian variety $H(3,q^2)$.  The point graph of $H(3,q^2)$ was considered in both~\cite{Gravier15} and~\cite{HegerGQ} with an upper bound of of $4q^2+q-5$ on the metric dimension.  This compares with our bound of $q^4+q^2+1$. 
\item For $[^2\! A_3(q)]$, we have the lines of $H(3,q^2)$, which is the dual of $Q^-(5,q)$.  The point graph of $Q^-(5,q)$ was considered in~\cite{Gravier15}, with an upper bound of $5q$ on the metric dimension.  In this case (with a change of variable from $\sqrt{q}$ to $q$) our bound evaluates to $q^4-q^2+1$. 
\item For $[^2\! A_4(q)]$, we have the lines of the Hermitian variety $H(4,q^2)$; while this is not the dual of a classical GQ, the dual was considered in~\cite{HegerGQ}, where an upper bound of $5q^3+1$ was obtained.  In comparison, our bound (again with a change of variable from $\sqrt{q}$ to $q$) is a polynomial of degree $7$.
\end{itemize}

The dual polar graphs on $[D_3(q)]$ are bipartite distance-regular graphs of diameter~$3$, so must be incidence graphs of symmetric designs; these turn out to be the point-hyperplane designs in $\mathrm{PG}(3,q)$.  General bounds for the metric dimension of incidence graphs of symmetric designs were obtained by the first author in~\cite[Corollary 2.5]{incidence}, while $\mathrm{PG}(3,q)$ (and, more generally, $\mathrm{PG}(n,q)$) were studied more closely by Bartoli {\em et al.}\ in~\cite{PGnq}.  In the latter paper, it is shown that there is a resolving set of size $8q$ in the incidence graph of $\mathrm{PG}(3,q)$; this improves on the general bounds for incidence graphs in~\cite{incidence}, which evaluate as $O(q\log q)$ for a design with these parameters, and are clearly much better than our bound of $O(q^3)$ for a graph with $O(q^3)$ vertices.\footnote{It would be interesting to see if an $O(q)$ bound holds for incidence graphs of arbitrary symmetric designs with the parameters of $\mathrm{PG}(3,q)$.}

If examined more closely, Theorem~\ref{theorem:main} is actually a result concerning the class dimension of association schemes arising from polar spaces. (See~\cite[Section 3.4]{bsmd} for more details on class dimension.)  Indeed, if minded so, given a polar space with parameters $q,d,e$ and some $t\in \{1,\ldots,d\}$, one might consider the graph having the collection of the totally isotropic subspaces of dimension $t$ as vertices, and where two such subspaces are adjacent if their intersection has codimension $1$. In particular, the dual polar graphs correspond to the case $t=d$. With minor modifications, the argument in the proof of Theorem~\ref{theorem:main} yields a bound on the $|\Omega_1|\times |\Omega_t|$ incidence matrix of the corresponding incidence structure. However, we have preferred to phrase our results only for the case $t=d$, i.e.\ for dual polar graphs, because only in this case the graph is distance-regular and hence only in this case our bound on the rank yields a bound on the metric dimension.

As far as we are aware, the problem of bounding from below the metric dimension of dual polar graphs $\Gamma(q,d,e)$ is entirely open.  Some relevant lower bounds in the case of generalized quadrangles (where $d=2$) are given in~\cite{Gravier15}.  It would be interesting to obtain a more general lower bound for the case of dual polar graphs of arbitrary diameter.

\subsection*{Acknowledgements}
The first author acknowledges financial support from an NSERC Discovery Grant and a Memorial University of Newfoundland startup grant.  The authors would like to thank John Bamberg and Maarten De Boeck for useful discussions.

\end{document}